# MULTIVARIATE NORMAL APPROXIMATION WITH STEIN'S METHOD OF EXCHANGEABLE PAIRS UNDER A GENERAL LINEARITY CONDITION

By Gesine Reinert[1] and Adrian Röllin[2]

*University of Oxford and National University of Singapore*

In this paper we establish a multivariate exchangeable pairs approach within the framework of Stein's method to assess distributional distances to potentially singular multivariate normal distributions. By extending the statistics into a higher-dimensional space, we also propose an embedding method which allows for a normal approximation even when the corresponding statistics of interest do not lend themselves easily to Stein's exchangeable pairs approach. To illustrate the method, we provide the examples of runs on the line as well as double-indexed permutation statistics.

**1. Introduction.** Stein's method was first published in Stein (1972) to assess the distance between univariate random variables and the normal distribution. This method has proved particularly powerful in the presence of both local dependence and weak global dependence.

A coupling at the heart of Stein's method for univariate normal approximation is the method of exchangeable pairs; see Stein (1986). Assume that $W$ is a univariate random variable with $\mathbb{E}W = 0$ and $\mathbb{E}W^2 = 1$, and assume that $W'$ is a random variable such that $(W, W')$ makes an exchangeable pair. Assume further that there is a number $\lambda > 0$ such that the conditional expectation of $W' - W$ with respect to $W$ satisfies

$$\mathbb{E}^W(W' - W) = -\lambda W. \tag{1.1}$$

Heuristically, (1.1) can be understood as a linear regression condition. If $(W, W')$ were bivariate normal with correlation $\rho$, then

$$\mathbb{E}^W W' = \rho W,$$

Received June 2008; revised December 2008.

[1]Supported in part by MMCOMNET Grant FP6-2003-BEST-Path-012999.

[2]Supported by Swiss National Science Foundation Project PBZH2–117033.

*AMS 2000 subject classifications.* Primary 60F05; secondary 62E17.

*Key words and phrases.* Stein's method, multivariate normal approximation, exchangeable pairs.

and (1.1) would be satisfied with $\lambda = 1 - \rho$. If $W$ was close to normal, then so would be $W'$, and it would not be unreasonable to assume that (1.1) is close to being satisfied.

In this spirit, the univariate theorem of Stein (1986) has been extended by Rinott and Rotar (1997). With the same basic setup as in Stein (1986), they generalize (1.1) by assuming that there is a number $\lambda > 0$ and a random variable $R = R(W)$ such that

$$\mathbb{E}^W(W' - W) = -\lambda W + R. \tag{1.2}$$

Note that, unlike condition (1.1), this is not a condition in the strict sense, as we can define $R := \mathbb{E}^W(W' - W) + \lambda W$ for any $\lambda$; however, we always have $\mathbb{E}R = 0$.

One of the results of Rinott and Rotar (1997) is that

$$\begin{aligned} &\sup_x |\mathbb{P}[W \le x] - \mathbb{P}[Z \le x]| \\ &\qquad \le \frac{6}{\lambda}\sqrt{\operatorname{Var}\mathbb{E}^W(W' - W)^2} + \frac{6}{\lambda^{1/2}}\sqrt{\mathbb{E}|W' - W|^3} + \frac{19}{\lambda}\sqrt{\operatorname{Var} R}, \end{aligned} \tag{1.3}$$

where $Z$ has standard normal distribution. So clearly, representation (1.2) is useful only if $\lambda^{-1}\sqrt{\operatorname{Var} R} = \mathrm{o}(1)$. In this case, if $\lambda_1$ and $\lambda_2$ stem from two different representations (1.2) for which $\lambda_i^{-1}\sqrt{\operatorname{Var} R_i} = \mathrm{o}(1)$ for $i = 1, 2$, then it it easy to see that $|\lambda_1 - \lambda_2|/(\lambda_1 + \lambda_2) = \mathrm{o}(1)$; in this sense, $\lambda$ is asymptotically unique. Rinott and Rotar (1997) then apply bound (1.3) to the number of ones in the anti-voter model, and to weighted $U$-statistics. Röllin (2008) provides a proof of a variant of (1.3) which does not use exchangeability but only $\mathscr{L}(W') = \mathscr{L}(W)$.

Stein's method has been extended to many other distributions; for an overview, see, for example, Reinert (2005). For multivariate normal approximations the method was first adapted by Barbour (1990) and Götze (1991), viewing the normal distribution as the stationary distribution of an Ornstein–Uhlenbeck diffusion, and using the generator of this diffusion as a characterizing operator for the normal distribution. Subsequent authors have used this generator approach for multivariate normal approximation with different variants, such as the local approach and the size-biasing approach by Goldstein and Rinott (1996) and Rinott and Rotar (1996), and the zero-biasing approach by Goldstein and Reinert (2005).

The exchangeable pair approach, in contrast, while having proved useful in non-normal contexts [see Chatterjee, Diaconis and Meckes (2005), Chatterjee, Fulman and Röllin (2006) and Röllin (2007)] remained restricted to the one-dimensional setting until very recently. A main stumbling block was that the extension of condition (1.2) to the multivariate setting is not obvious from the viewpoint of Stein's method.

In Chatterjee and Meckes (2008), this issue was finally addressed. They propose the condition that, for all $i = 1, \dots, d$,

$$\mathbb{E}^W(W_i' - W_i) = -\lambda W_i, \tag{1.4}$$

for a fixed number $\lambda$, where now $W = (W_1, \dots, W_d)$ and $W' = (W_1', \dots, W_d')$ are identically distributed $d$-vectors with uncorrelated components (an extension to the additional remainder term $R$ was not considered, but would be straightforward). They employ such couplings to bound the distance to the standard multivariate normal distribution. Using the same argument as Röllin (2008), Chatterjee and Meckes (2008) are able to give proofs of their theorems without using exchangeability and apply them successfully to various multivariate applications.

Applying a similar heuristic as for (1.1), however, if $(W, W')$ were jointly normal, with mean vector 0 and covariance matrix

$$\Sigma_0 = \begin{pmatrix} \Sigma & \tilde{\Sigma} \\ \tilde{\Sigma} & \Sigma \end{pmatrix}, \tag{1.5}$$

then $\mathbb{E}^W W' = \tilde{\Sigma}\Sigma^{-1}W$ [see, e.g., Mardia, Kent and Bibby (1979), page 63, Theorem 3.2.4], in which case

$$\mathbb{E}^W(W' - W) = -(\mathrm{Id} - \tilde{\Sigma}\Sigma^{-1})W; \tag{1.6}$$

here Id denotes the identity matrix. Again, if $(W, W')$ is approximately jointly normal, then we expect (1.6) to be approximately satisfied. This heuristic leads to the condition that

$$\mathbb{E}^W(W' - W) = -\Lambda W + R \tag{1.7}$$

for an invertible $d \times d$ matrix $\Lambda$ and a remainder term $R = R(W)$. For $R = 0$, even if $\Sigma = \mathrm{Id}$, we would obtain $\Lambda = \mathrm{Id} - \tilde{\Sigma}$, which in general is not diagonal. Hence, we argue that (1.7) is not only more general, but also more natural than (1.4).

Different exchangeable pairs will lead to different $\Lambda$ and $R$ in (1.7); our embedding method suggests suitable decompositions. Indeed, for a specific exchangeable pair $(W, W')$ at hand, it is often far from obvious whether this pair will satisfy the linearity condition (1.7) with $R$ of the required small order, unless equal to zero. Consider the case of 2-runs. For a sequence of i.i.d. Bernoulli distributed random variables $\xi_1, \dots, \xi_n$ such that $\mathbb{P}[\xi_1 = 1] = p$, define the centered number of 2-runs

$$V_2 = \sum_{i=1}^{n} \xi_i \xi_{i+1} - np^2,$$

where we let $\xi_{n+1} := \xi_1$. The most natural construction of an exchangeable pair in the spirit of Stein (1986) is to pick uniformly a $\xi_i$ and replace it by an independent copy $\xi_i'$. Denote by $V_2'$ the resulting number of 2-runs in the new sequence. It is easy to calculate (see Section 4.2) that

$$\mathbb{E}^{V_2}(V_2' - V_2) = -\frac{2}{n}V_2 + \frac{2p}{n}\mathbb{E}^{V_2}\sum_{i=1}^{n}(\xi_i - p). \tag{1.8}$$

The conditional expectation on the right-hand side of (1.8) is hard to calculate. Furthermore, it has the same order of magnitude as $V_2$. Also, the weighted $U$-statistics approach of Rinott and Rotar (1997) (Proposition 1.2) does not yield convergent bounds to the normal distribution. We propose the following approach to this problem. Keeping the above coupling, we define $V_1 := \sum_{i=1}^n \xi_i - np$ (and $V_1'$ accordingly) and consider the problem as a 2-dimensional problem $W := \binom{V_1}{V_2}$. Equation (1.8) now yields $\mathbb{E}^W(V_2' - V_2) = -\frac{2}{n}V_2 + \frac{2p}{n}V_1$, and further calculations reveal that $\mathbb{E}^W(V_1' - V_1) = -\frac{1}{n}V_1$, so that now (1.7) holds with

$$\Lambda = \frac{1}{n}\begin{bmatrix} 1 & 0 \\ -2p & 2 \end{bmatrix}$$

and $R = 0$. Using this embedding into a higher-dimensional setting, the problem now fits into our framework and allows not only for a normal approximation of the primary statistic, but for an approximation of the joint distribution of the primary and auxiliary statistics. For this embedding method, the generality of condition (1.7) is essential; see (4.1) later.

The rest of the article is organized as follows. In the next section we prove an abstract nonsingular multivariate normal approximation theorem (Theorem 2.1) for smooth test functions. The explicit bound on the distance to the normal distribution is given in terms of the conditional variance, the absolute third moments and the variance of the remainder term. Proposition 2.8 gives the extension to singular multivariate normal distributions, using Stein's method and the triangle inequality. To illustrate our results, we calculate the example of sums of i.i.d. variables.

Section 3 uses the abstract theorem to obtain a similar result for nonsmooth test functions, such as indicators of convex sets. Adapting the approach by Rinott and Rotar (1996) to general multivariate normal approximation, Corollary 3.1 displays how the main terms involved in the error bounds for smooth test functions reappear in the bounds for nonsmooth test functions.

Section 4 discusses the above mentioned embedding method and illustrates its application with a detailed treatment of runs on the line. We also sketch the application to double-indexed permutation statistics.

The generality of (1.7) comes at the extra cost that now exchangeability seems almost inevitable. Indeed, in view of Röllin (2008), we were surprised that, in the multivariate setting, the exchangeability condition cannot be removed as easily as in the one-dimensional case. Therefore, the last section discusses the exchangeability condition, condition (1.7) and their implications.

Appendix A contains the proof of Corollary 3.1, and details of the runs example are in Appendix B.

1.1. *Notation.* Random vectors in $\mathbb{R}^d$ are written in the form $W = (W_1, W_2, \dots, W_d)^t$, where $W_i$ are $\mathbb{R}$-valued random variables for $i = 1, \dots, d$. If $\Sigma$ is a symmetric, nonnegative definite matrix, we denote by $\Sigma^{1/2}$ the unique symmetric, nonnegative definite square root of $\Sigma$. Denote by Id the identity matrix, usually of dimension $d$. Throughout this article, $Z$ will denote a random vector having standard multivariate normal distribution, also of dimension $d$.

For ease of presentation, we abbreviate the transpose of the inverse of a matrix in the form $\Lambda^{-t} := (\Lambda^{-1})^t$.

Stein's method makes good use of Taylor expansions. For derivatives of smooth functions $h : \mathbb{R}^d \to \mathbb{R}$, we use the notation $\nabla$ for the gradient operator. For the sake of presentation, the partial derivatives are abbreviated as $h_i = \frac{\partial}{\partial x_i} h$, $h_{i,j} = \frac{\partial^2}{\partial x_i \partial x_j} h$ unless we would like to emphasise the dependence on the variables.

To derive uniform bounds, we shall employ the supremum norm, denoted by $\|\cdot\|$ for both functions and matrices. For a function $h : \mathbb{R}^d \to \mathbb{R}$, we abbreviate $|h|_1 := \sup_i \|\frac{\partial}{\partial x_i} h\|$, $|h|_2 := \sup_{i,j} \|\frac{\partial^2}{\partial x_i \partial x_j} h\|$, and so on, if the corresponding derivatives exist.

**2. The distance to multivariate normal distribution in terms of smooth test functions.** First we derive a bound on the distance between a multivariate target distribution and a multivariate normal distribution with the same mean vector (which is assumed to be 0 in the sequel), and with the same, positive definite covariance matrix. We start by considering smooth test functions.

THEOREM 2.1. *Assume that $(W, W')$ is an exchangeable pair of $\mathbb{R}^d$-valued random vectors such that*

$$\mathbb{E}W = 0, \qquad \mathbb{E}WW^t = \Sigma, \tag{2.1}$$

*with $\Sigma \in \mathbb{R}^{d\times d}$ symmetric and positive definite. Suppose further that (1.7) is satisfied for an invertible matrix $\Lambda$ and a $\sigma(W)$-measurable random vector $R$. Then, if $Z$ has $d$-dimensional standard normal distribution, we have for every three times differentiable function $h$,*

$$|\mathbb{E}h(W) - \mathbb{E}h(\Sigma^{1/2}Z)| \le \frac{|h|_2}{4} A + \frac{|h|_3}{12} B + \left(|h|_1 + \frac{1}{2} d \|\Sigma\|^{1/2} |h|_2\right) C, \tag{2.2}$$

*where, with* $\lambda^{(i)} := \sum_{m=1}^d |(\Lambda^{-1})_{m,i}|$,

$$A = \sum_{i,j=1}^d \lambda^{(i)} \sqrt{\operatorname{Var} \mathbb{E}^W (W_i' - W_i)(W_j' - W_j)},$$

$$B = \sum_{i,j,k=1}^d \lambda^{(i)} \mathbb{E}|(W_i' - W_i)(W_j' - W_j)(W_k' - W_k)|,$$

$$C = \sum_{i=1}^d \lambda^{(i)} \sqrt{\operatorname{Var} R_i}.$$

Before we proceed with the proof, we illustrate Theorem 2.1 by means of the simple example of sums of i.i.d. random variables and make some further remarks.

COROLLARY 2.2. *Suppose that* $W = (W_1, \ldots, W_d)$ *is such that, for each* $i$, $W_i = \sum_{j=1}^n X_{i,j}$, *where* $X_{i,j}, i = 1, \ldots, d, j = 1, \ldots, n$, *are i.i.d. with mean zero and variance* $\frac{1}{n}$, *so that the covariance matrix* $\Sigma = \mathrm{Id}$. *Assume further that there exist* $0 < \beta, \gamma < \infty$ *such that*

$$\mathbb{E}|X_{i,j}|^3 = \beta n^{-3/2} \quad \text{and} \quad \operatorname{Var}(X_{i,j}^2) = \gamma n^{-2}.$$

*Then, for every three times differentiable function* $h$,

$$|\mathbb{E}h(W) - \mathbb{E}h(Z)| \le \frac{d}{\sqrt{n}} \left( \frac{\sqrt{\gamma}}{4} |h|_2 + \frac{2\beta}{3} |h|_3 \right).$$

PROOF. We construct an exchangeable pair by choosing a vector $I$ and a summand $J$ uniformly, such that $\mathbb{P}(I = i, J = j) = 1/(dn)$. If $I = i, J = j$, we replace $X_{i,j}$ by an independent copy $X'_{i,j}$; all other variables remain unchanged. Put

$$W_I' = W_I - X_{I,J} + X'_{I,J}$$

and $W_k' = W_k$ for $k \ne I$; denote by $W'$ the resulting $d$-vector. Then $(W, W')$ is exchangeable, and, in (1.7), $\Lambda = \frac{1}{dn} \mathrm{Id}$ with $R = 0$ and, hence, $C = 0$. For our bounds we note that $\lambda^{(i)} = dn$. We calculate that

$$\mathbb{E}^W (W_i' - W_i)^2 = \frac{1}{dn} + \frac{1}{dn} \sum_j \mathbb{E}^W X_{i,j}^2.$$

Thus,

$$\operatorname{Var} \mathbb{E}^W (W_i' - W_i)^2 \le \frac{1}{d^2 n^2} \sum_j \operatorname{Var} X_{ij}^2 \le \frac{\gamma}{n^3 d^2}.$$

Moreover, by the construction, for $i \neq k$, $(W_i' - W_i)(W_k' - W_k) = 0$, and $(W_i' - W_i)(W_k' - W_k)(W_l' - W_l) = 0$, unless $i = k = l$. By assumption,

$$\mathbb{E}|W_i' - W_i|^3 = \frac{1}{dn} \sum_{\ell=1}^{d} \mathbf{1}(\ell = i) \sum_{j=1}^{n} \mathbb{E}|X_{i,j} - X_{i,j}'|^3 \leq \frac{8\beta}{dn^{3/2}}.$$

The result now follows directly from Theorem 2.1. □

Remark 2.3. Multivariate normal approximations for vectors of sums of i.i.d. random variables have been so intensively studied that there is not enough space to review all the results. The approach most similar to ours is found in Chatterjee and Meckes (2008), where instead of exchanging only one summand, a whole vector would be exchanged. Their results yield

$$|\mathbb{E}h(W) - \Phi h| \leq \frac{d^{3/2}\sqrt{\gamma + 1}}{2\sqrt{n}} |h|_1 + 4\frac{d^3\beta}{\sqrt{n}} |h|_2.$$

Due to the different Stein equation used, the dependence on the dimension differs, and the bounds are in terms of different derivatives of the test function. The overall similarity in this special case is apparent.

Remark 2.4. If we were to normalize the random vectors in Theorem 2.1, denoting the normalization of $W$ by $\hat{W} := \Sigma^{-1/2} W$ and $\hat{W}' = \Sigma^{-1/2} W'$, then, the conditions of the theorem remain satisfied for $(\hat{W}, \hat{W}')$ with $\hat{\Sigma} = \mathrm{Id}$ and $\hat{\Lambda} = \Sigma^{-1/2} \Lambda \Sigma^{1/2}$ as well as $\hat{R} = \Sigma^{-1/2} R$.

Remark 2.5. As a precursor to (1.7), in the context of multivariate zero-biasing, Goldstein and Reinert (2005) use the condition of the form (1.7) for $\Lambda$ such that $\Lambda_{i,j} = \rho + \mathbf{1}(i = j)$.

After these remarks we proceed to the proof of Theorem 2.1, which is based on the Stein characterization of the normal distribution that $Y \in \mathbb{R}^d$ is a multivariate normal $\mathrm{MVN}(0, \Sigma)$ if and only if

$$\mathbb{E}\{\nabla^t \Sigma \nabla f(Y) - Y^t \nabla f(Y)\} = 0 \qquad \text{for all smooth} f : \mathbb{R}^d \to \mathbb{R}. \tag{2.3}$$

We will need the following lemma to prove the theorem; however, see also Remark 2.4, Barbour (1990), Goldstein and Rinott (1996) and Götze (1991). The proof of Lemma 2.6 is routine.

Lemma 2.6. *Let $h : \mathbb{R}^d \to \mathbb{R}$ be differentiable with bounded first derivative. Then, if $\Sigma \in \mathbb{R}^{d \times d}$ is symmetric and positive definite, there is a solution $f : \mathbb{R}^d \to \mathbb{R}$ to the equation*

$$\nabla^t \Sigma \nabla f(w) - w^t \nabla f(w) = h(w) - \mathbb{E}h(\Sigma^{1/2} Z), \tag{2.4}$$

*which holds for every $w \in \mathbb{R}^d$. If, in addition, $h$ is $n$ times differentiable, there is a solution $f$ which is also $n$ times differentiable and we have for every $k = 1, \dots, n$ the bound*

$$\left| \frac{\partial^k f(w)}{\prod_{j=1}^k \partial w_{i_j}} \right| \leq \frac{1}{k} \left| \frac{\partial^k h(w)}{\prod_{j=1}^k \partial w_{i_j}} \right| \tag{2.5}$$

*for every $w \in \mathbb{R}^d$.*

REMARK 2.7. Compared to the main theorem of Chatterjee and Meckes (2008), which only needs the existence of two derivatives, our Theorem 2.1 is more restrictive in the choice of test functions $h$. This reflects the fact that we make use of Lemma 2.6, which is motivated by Goldstein and Rinott (1996), whereas Chatterjee and Meckes (2008) prove new bounds on the solutions of (2.4), but only for $\Sigma = \mathrm{Id}$; see also Raič (2004) for similar results. The general result of Lemma 2.6, however, allows us to work with the unstandardized pair $(W, W')$, which not only usually simplifies the calculations, but also yields more informative bounds if the limiting covariance matrix is singular.

PROOF OF THEOREM 2.1. Our aim is to bound $|\mathbb{E}h(W) - \mathbb{E}h(\Sigma^{1/2}Z)|$ by bounding $|\mathbb{E}\{\nabla^t \Sigma \nabla f(W) - W^t \nabla f(W)\}|$, where $f$ is the solution to the Stein equation (2.4). First we expand $\mathbb{E}W^t \nabla f(W)$. Define the real-valued, anti-symmetric function

$$F(w', w) := \tfrac{1}{2}(w' - w)^t \Lambda^{-t}(\nabla f(w') + \nabla f(w)) \tag{2.6}$$

for $w, w' \in \mathbb{R}^d$, and note that, because of exchangeability, $\mathbb{E}F(W', W) = 0$; see Stein (1986). Thus,

$$\begin{aligned} 0 &= \tfrac{1}{2}\mathbb{E}\{(W' - W)^t \Lambda^{-t}(\nabla f(W') + \nabla f(W))\} \\ &= \mathbb{E}\{(W' - W)^t \Lambda^{-t} \nabla f(W)\} \\ &\quad + \tfrac{1}{2}\mathbb{E}\{(W' - W)^t \Lambda^{-t}(\nabla f(W') - \nabla f(W))\} \\ &= \mathbb{E}\{R^t \Lambda^{-t} \nabla f(W)\} - \mathbb{E}\{W^t \nabla f(W)\} \\ &\quad + \tfrac{1}{2}\mathbb{E}\{(W' - W)^t \Lambda^{-t}(\nabla f(W') - \nabla f(W))\}, \end{aligned} \tag{2.7}$$

where we used (1.7) for the last step. Recalling the notation $f_{i,j}(x) = \frac{\partial^2}{\partial x_i \partial x_j} f(x)$, Taylor expansion gives

$$\begin{aligned} &(w' - w)^t \Lambda^{-t}(\nabla f(w') - \nabla f(w)) \\ &\quad = \sum_{m,i,j} (\Lambda^{-1})_{m,i}(w'_i - w_i)(w'_j - w_j) f_{m,j}(w) \\ &\qquad + \sum_{m,i,j,k} (\Lambda^{-1})_{m,i}(w'_i - w_i)(w'_j - w_j)(w'_k - w_k)\tilde{R}_{mjk}, \end{aligned}$$

where

$$|\tilde{R}_{mjk}| \le \frac{1}{2}\left\|\frac{\partial^3 f}{\partial w_m\,\partial w_j\,\partial w_k}\right\|. \tag{2.8}$$

Thus, in (2.7),

$$\begin{aligned}
&\mathbb{E}\{(W'-W)^t\Lambda^{-t}(\nabla f(W')-\nabla f(W))\}\\
&\quad= \sum_{m,i,j}(\Lambda^{-1})_{m,i}\mathbb{E}\{(W_i'-W_i)(W_j'-W_j)f_{m,j}(W)\}\\
&\qquad+\sum_{m,i,j,k}(\Lambda^{-1})_{m,i}\mathbb{E}\{(W_i'-W_i)(W_j'-W_j)(W_k'-W_k)\tilde{R}_{mjk}\}.
\end{aligned}\tag{2.9}$$

Now we turn our attention to $\mathbb{E}\nabla^t\Sigma\nabla f(W)$. Note that, because of (2.1), (1.7) and exchangeability,

$$\begin{aligned}
&\mathbb{E}(W'-W)(W'-W)^t\\
&\quad= \mathbb{E}\{W(W-W')^t\}+\mathbb{E}\{W(W-W')^t\}\\
&\quad= 2\mathbb{E}\{W(\Lambda W-R)^t\}=2\Sigma\Lambda^t-2\mathbb{E}(WR^t)=:T.
\end{aligned}\tag{2.10}$$

Hence,

$$\begin{aligned}
\nabla^t\Sigma\nabla f(w) &= \frac{1}{2}\nabla^t T\Lambda^{-t}\nabla f(w)+\nabla^t\mathbb{E}(WR^t)\Lambda^{-t}\nabla f(w)\\
&= \frac{1}{2}\sum_{m,i,j}(\Lambda^{-1})_{m,i}T_{j,i}\frac{\partial^2 f(w)}{\partial w_m\,\partial w_j}+\sum_{m,i,j}(\Lambda^{-1})_{m,i}\mathbb{E}(W_jR_i)\frac{\partial^2 f(w)}{\partial w_m\,\partial w_j}.
\end{aligned}$$

Combining this equation with (2.7) and (2.9),

$$\begin{aligned}
&|\mathbb{E}\{\nabla^t\Sigma\nabla f(W)-W^t\nabla f(W)\}|\\
&\quad\le \frac{1}{2}\left|\sum_{m,i,j}\mathbb{E}\left\{(\Lambda^{-1})_{m,i}[T_{j,i}-\mathbb{E}^W(W_i'-W_i)(W_j'-W_j)]\frac{\partial^2 f(W)}{\partial w_m\,\partial w_j}\right\}\right|\\
&\qquad+\frac{1}{2}\left|\sum_{m,i,j,k}\mathbb{E}\{(\Lambda^{-1})_{m,i}(W_i'-W_i)(W_j'-W_j)(W_k'-W_k)\tilde{R}_{mjk}\}\right|\\
&\qquad+\left|\sum_{i,m}(\Lambda^{-1})_{m,i}\mathbb{E}\left\{R_i\frac{\partial f(W)}{\partial w_m}\right\}\right|+\left|\sum_{m,i,j}(\Lambda^{-1})_{m,i}\mathbb{E}(W_jR_i)\mathbb{E}\left\{\frac{\partial^2 f(W)}{\partial w_m\,\partial w_j}\right\}\right|\\
&\quad\le \frac{|h|_2}{4}\sum_{i,j}\lambda^{(i)}\mathbb{E}|T_{j,i}-\mathbb{E}^W(W_i'-W_i)(W_j'-W_j)|+\frac{|h|_3}{12}B\\
&\qquad+|h|_1\sum_i\lambda^{(i)}\mathbb{E}|R_i|+\frac{|h|_2}{2}\sum_{i,j}\lambda^{(i)}\mathbb{E}|W_jR_i|,
\end{aligned}\tag{2.11}$$

where we used (2.8) to obtain the second inequality, and Lemma 2.6 to obtain the last inequality. From the Cauchy–Schwarz inequality, $\mathbb{E}|R_j| \le \sqrt{\mathbb{E}R_j^2}$ and

$$\mathbb{E}|W_j R_i| \le \sqrt{\mathbb{E}W_j^2 \mathbb{E}R_i^2} \le \|\Sigma\|^{1/2}\sqrt{\mathbb{E}R_i^2}.$$

The $C$-expression in (2.2) now follows from the last two terms of (2.11). Recalling that $\mathbb{E}(W' - W)(W' - W)^t = T$, this proves the first term of (2.2) from the first term of (2.11). □

Sometimes we may wish to assess the distance to a normal distribution for which the covariance matrix $\Sigma_0$, while nonnegative definite, does not have full rank. Stein's method helps to derive a straightforward bound in this case also. The proof of the following proposition is straightforward and routine, noting that (2.3) remains valid if the covariance matrix is not of full rank.

PROPOSITION 2.8. *Let $X$ and $Y$ be $\mathbb{R}^d$-valued normal vectors with distributions $X \sim \mathrm{MVN}(0, \Sigma)$ and $Y \sim \mathrm{MVN}(0, \Sigma_0)$, where $\Sigma = (\sigma_{i,j})_{i,j=1,\dots,d}$ has full rank, and $\Sigma_0 = (\sigma^0_{i,j})_{i,j=1,\dots,d}$ is nonnegative definite. Let $h : \mathbb{R}^d \to \mathbb{R}$ have 2 bounded derivatives. Then*

$$|\mathbb{E}h(X) - \mathbb{E}h(Y)| \le \frac{1}{2}|h|_2 \sum_{i,j=1}^{d} |\sigma_{i,j} - \sigma^0_{i,j}|.$$

Using the triangle inequality and Theorem 2.1, we thus obtain a bound for a normal approximation even for a normal distribution with degenerate covariance matrix.

**3. Nonsmooth test functions.** Following Rinott and Rotar (1996), let $\Phi$ denote the standard normal distribution in $\mathbb{R}^d$, and $\phi$ the corresponding density function. For $h : \mathbb{R}^d \to \mathbb{R}$ set

$$h^+_\delta(x) = \sup\{h(x+y) : |y| \le \delta\},$$
$$h^-_\delta(x) = \inf\{h(x+y) : |y| \le \delta\},$$
$$\tilde{h}(x, \delta) = h^+_\delta(x) - h^-_\delta(x).$$

Let $\mathcal{H}$ be a class of measurable functions $\mathbb{R}^d \to \mathbb{R}$ which are uniformly bounded by 1. Suppose that, for any $h \in \mathcal{H}$:

(C1) for any $\delta > 0$, $h^+_\delta(x)$ and $h^-_\delta(x)$ are in $\mathcal{H}$,
(C2) for any $d \times d$ matrix $A$ and any vector $b \in \mathbb{R}^d$, $h(Ax + b) \in \mathcal{H}$,

(C3) for some constant $a = a(\mathcal{H}, \delta)$

$$\sup_{h\in\mathcal{H}}\left\{\int_{\mathbb{R}^d}\tilde{h}(x,\delta)\Phi(dx)\right\} \le a\delta. \tag{3.1}$$

Obviously we may assume $a \ge 1$.

The class of indicators of measurable convex sets is such a class; for this class, $a \le 2\sqrt{d}$; see Bolthausen and Götze (1993).

In the same way as in Rinott and Rotar (1996), we can show the following corollary. The presentation differs from Rinott and Rotar (1996), as we make the relationship to the bounds in Theorem 2.1 immediate and in that we allow for general $\Sigma$. The now fairly standard proof is found in Appendix A. We also note forthcoming work by Bhattacharya and Holmes (2007).

Let $W$ have mean vector 0 and variance–covariance matrix $\Sigma$. If $\Lambda$ and $R$ are such that (1.7) is satisfied for $W$, then $Y = \Sigma^{-1/2}W$ satisfies (1.7) with $\hat{\Lambda} = \Sigma^{-1/2}\Lambda\Sigma^{1/2}$ and $R' = \Sigma^{-1/2}R$. We put

$$\hat{\lambda}^{(i)} = \sum_{m=1}^{d}|(\Sigma^{-1/2}\Lambda^{-1}\Sigma^{1/2})_{m,i}|,$$

as well as

$$A' = \sum_{i,j}\hat{\lambda}^{(i)}\sqrt{\operatorname{Var}\mathbb{E}^Y\sum_{k,\ell}\Sigma_{i,k}^{-1/2}\Sigma_{j,\ell}^{-1/2}(W_k'-W_k)(W_\ell'-W_\ell)},$$

$$B' = \sum_{i,j,k}\hat{\lambda}^{(i)}\mathbb{E}\left|\sum_{r,s,t}\Sigma_{i,r}^{-1/2}\Sigma_{j,s}^{-1/2}\Sigma_{k,t}^{-1/2}(W_r'-W_r)(W_s'-W_s)(W_t'-W_t)\right|$$

and

$$C' = \sum_{i=1}^{d}\hat{\lambda}^{(i)}\sqrt{\mathbb{E}\left(\sum_k\Sigma_{i,k}^{-1/2}R_k\right)^2}.$$

Corollary 3.1. *Let $W$ be as in Theorem 2.1. Then, for all $h \in \mathcal{H}$ with $|h| \le 1$, there exists $\gamma = \gamma(d)$ such that, with $a > 1$ as in (3.1),*

$$\sup_{h\in\mathcal{H}}|\mathbb{E}h(W) - \mathbb{E}h(Z)| \le \gamma^2\left(-D'\log(T') + \frac{B'}{2\sqrt{T'}} + C' + a\sqrt{T'}\right),$$

*with*

$$T' = \frac{1}{a^2}\left(D' + \sqrt{\frac{aB'}{2} + D'^2}\right)^2 \quad \text{and} \quad D' = \frac{A'}{2} + C'd.$$

If $A', B'$ and $C'$ are $O(n^{-1/2})$, then we would obtain a bound of order $O(n^{-1/4})$. This is poorer than the $n^{-1/2}\log n$ type of bounds obtained in Rinott and Rotar (1996), but Rinott and Rotar (1996) obtain the improved rate by assuming that the random vectors are bounded.

## 4. The embedding method and applications.

4.1. *General framework.* Assume that an $\ell$-dimensional random vector $W_{(\ell)}$ of interest is given. Often, the construction of an exchangeable pair $(W_{(\ell)}, W'_{(\ell)})$ is straightforward. If, say, $W_{(\ell)} = W_{(\ell)}(\mathbb{X})$ is a function of i.i.d. random variables $\mathbb{X} = (X_1, \ldots, X_n)$, one can choose uniformly an index $I$ from 1 to $n$, replace $X_I$ by an independent copy $X'_I$, and define $W'_{(\ell)} := W_{(\ell)}(\mathbb{X}')$, where $\mathbb{X}'$ is now the vector $\mathbb{X}$ but with $X_I$ replaced by $X'_I$.

In general there is no hope that $(W_{(\ell)}, W'_{(\ell)})$ will satisfy condition (1.2) with $R$ being of the required smaller order or even equal to zero, so that in this case Theorem 2.1 would not yield useful bounds.

Surprisingly often it is possible, though, to extend $W_{(\ell)}$ to a vector $W \in \mathbb{R}^d$ such that we can construct an exchangeable pair $(W, W')$ which satisfies condition (1.2) with $R = 0$. If we can bound the distance of the distribution $\mathcal{L}(W)$ to a $d$-dimensional multivariate normal distribution, then a bound on the distance of the distribution $\mathcal{L}(W_{(\ell)})$ to an $\ell$-dimensional multivariate normal distribution follows immediately.

To explain the approach, we turn the problem on its head. Suppose that $W \in \mathbb{R}^d$ is such that we can construct an exchangeable pair $(W, W')$ which satisfies condition (1.2) with $R = 0$. Rename the first $\ell$ components to comprise $W_{(\ell)}$, so that

$$W = \begin{bmatrix} W_{(\ell)} \\ W^{(d-\ell)} \end{bmatrix},$$

and $W_{(\ell)} = I_{\ell,0} W$, with

$$I_{\ell,0} = (\mathrm{Id}_\ell, 0_{\ell\times(d-\ell)}),$$

$0_{\ell\times(d-\ell)}$ denoting the $\ell \times (d-\ell)$-matrix consisting entirely of 0's. Defining $W'_{(\ell)} = I_{\ell,0} W'$, it follows that $(W_{(\ell)}, W'_{(\ell)})$ forms an exchangeable pair. From (1.2),

$$\mathbb{E}^W (W_{(\ell)} - W'_{(\ell)}) = I_{\ell,0} \mathbb{E}^W (W - W') = -I_{\ell,0} \Lambda W.$$

Now decompose the matrix $\Lambda$ as

$$\Lambda = \begin{bmatrix} \Lambda_{1,1} & \Lambda_{1,2} \\ \Lambda_{2,1} & \Lambda_{2,2} \end{bmatrix},$$

where $\Lambda_{1,1}$ denotes an $\ell \times \ell$ submatrix, $\Lambda_{1,2}$ denotes an $\ell \times (d-\ell)$ submatrix, and so on. Then

$$I_{\ell,0} \Lambda W = \Lambda_{1,1} W_{(\ell)} + \Lambda_{1,2} W^{(d-\ell)}$$

and, hence,

$$\mathbb{E}^W (W_{(\ell)} - W'_{(\ell)}) = -\Lambda_{1,1} W_{(\ell)} - \Lambda_{1,2} W^{(d-\ell)}.$$

Conditioning on $W_{(\ell)}$ gives that

$$\mathbb{E}^{W_{(\ell)}}(W_{(\ell)} - W'_{(\ell)}) = -\Lambda_{1,1}W_{(\ell)} - \Lambda_{1,2}\mathbb{E}^{W_{(\ell)}}W^{(d-\ell)}.$$

Thus, condition (1.2) is satisfied with

$$R = -\Lambda_{1,2}\mathbb{E}^{W_{(\ell)}}W^{(d-\ell)}. \tag{4.1}$$

If $\Lambda_{1,2} = 0$, then no embedding is required. But if $\Lambda_{1,2} \neq 0$, then the remainder $R$ in (1.2) is a nontrivial linear combination of random variables, and these random variables could serve as embedding vector. In order to obtain useful bounds in Theorem 2.1, the embedding dimension $d$ should not be too large. In the examples below it will be obvious how to choose $W^{(d-\ell)}$ to make the construction work.

4.2. *Runs on the line.* Let $\mathbb{X} = (\xi_1, \dots, \xi_n)$ be a sequence of independent random variables with distribution Bernoulli$(p)$, $0 < p < 1$, that is, $\mathbb{P}[\xi_1 = 1] = 1 - \mathbb{P}[\xi_1 = 0] = p$. For $d > 1$, define the (centered) number of $d$-runs as

$$V_d := \sum_{m=1}^{n} (\xi_m \xi_{m+1} \cdots \xi_{m+d-1} - p^d),$$

where, for convenience, we assume the torus convention that $\xi_{n+1} \equiv \xi_1$, $\xi_{n+2} \equiv \xi_2$ and so on.

As mentioned in the Introduction, if we want to use the obvious construction of an exchangeable pair, the univariate version of exchangeable pairs of Rinott and Rotar (1997) (Proposition 1.2) does not yield convergent bounds of $V_d$ to the standard normal distribution if $d > 1$. However, we can tackle the example with our approach by incorporating the auxiliary variables $V_1, \dots, V_{d-1}$, such that the problem becomes linear in a higher-dimensional setting.

We construct an exchangeable pair $(\mathbb{X}, \mathbb{X}')$, where instead of just one, we resample $d-1$ of the $\xi_i$. To this end, let $I$ be uniformly distributed over $\{1, \dots, n\}$ and let $\tilde{\xi}_1, \dots, \tilde{\xi}_n$ be independent copies of the $\xi_i$. Let $\mathbb{X}'$ be the same as $\mathbb{X}$ but with the subsequence $\xi_I, \xi_{I+1}, \dots, \xi_{I+d-2}$ of length $d-1$ replaced by $\xi'_I, \xi'_{I+1}, \dots, \xi'_{I+d-2}$. Clearly $(\mathbb{X}, \mathbb{X}')$ forms an exchangeable pair. Define $V'_i := V_i(\mathbb{X}')$; we have

$$\begin{aligned} V'_i - V_i = {} & \sum_{m=I-i+1}^{I-1} \xi_m \cdots \xi_{I-1}\xi'_I \cdots \xi'_{m+i-1} + \sum_{m=I}^{I+d-i-1} \xi'_m \cdots \xi'_{m+i-1} \\ & + \sum_{m=I+d-i}^{I+d-2} \xi'_m \cdots \xi'_{I-1}\xi_I \cdots \xi_{m+i-1} - \sum_{m=I-i+1}^{I+d-2} \xi_m \cdots \xi_{m+i-1}, \end{aligned} \tag{4.2}$$

where sums $\sum_a^b$ are defined to be zero if $a > b$. Now, (4.2) yields

$$\begin{aligned}
&\mathbb{E}^{(V_1,\ldots,V_{d-1})}(V_i' - V_i) \\
&\quad = -n^{-1}[(d+i-2)V_i - 2pV_{i-1} - 2p^2V_{i-2} - \cdots - 2p^{i-1}V_1] \\
&\quad = -n^{-1}\left[(d+i-2)V_i - 2\sum_{k=1}^{i-1} p^{i-k}V_k\right].
\end{aligned} \tag{4.3}$$

From this representation we see that we may take $V_1, \ldots, V_{d-1}$ as the auxiliary random variables.

Straightforward calculations yield that, for all $i \geq j$,

$$\begin{aligned}
\mathbb{E}(V_iV_j) &= n\left[(i-j+1)p^i + 2\sum_{l=1}^{j-1} p^{i+j-l} - (i+j-1)p^{i+j}\right] \\
&= np^i(1-p)\sum_{k=0}^{j-1}(i-j+1+2k)p^k.
\end{aligned} \tag{4.4}$$

In particular,

$$\mathbb{E}V_i^2 = np^i(1-p)\sum_{k=0}^{i-1}(1+2k)p^k, \tag{4.5}$$

which lies in the interval between $np^i(1-p)$ and $np^i(1-p)i^2$. Thus, we define the $W_i$ to be the weighted versions

$$W_i := \frac{V_i}{\sqrt{np^i(1-p)}}, \tag{4.6}$$

and from (4.4) we have for general $i$ and $j$

$$\mathbb{E}(W_iW_j) = p^{|i-j|/2}\sum_{k=0}^{i\wedge j-1}(|i-j|+1+2k)p^k =: \sigma_{i,j}. \tag{4.7}$$

From (4.7) it is clear that the corresponding $\Sigma = (\sigma_{i,j})_{i,j}$ is constant for all $n$ and of full rank. For $p \to 0$, $\Sigma$ converges to uncorrelated coordinates and for $p \to 1$ to a matrix of rank 1. For applications and further references see Glaz, Naus and Wallenstein (2001) and Balakrishnan and Koutras (2002). Now, from (4.3) we have

$$\mathbb{E}^W(W_i' - W_i) = -n^{-1}\left[(d+i-2)W_i - 2\sum_{k=1}^{i-1} p^{(i-k)/2}W_k\right].$$

Thus, (1.7) is satisfied with $R=0$ and

$$\Lambda = \frac{1}{n}\begin{bmatrix} d-1 & & & & & \\ -2p^{1/2} & d & & & 0 & \\ \vdots & & \ddots & & & \\ -2p^{(k-1)/2} & \cdots & -2p^{1/2} & d+k-2 & & \\ \vdots & & & & \ddots & \\ -2p^{(d-1)/2} & & \cdots & & -2p^{1/2} & 2(d-1) \end{bmatrix}.$$

THEOREM 4.1. *With $W$ defined as in (4.6), $n > 2d-1$ and $\Sigma$ given through (4.7), we have for three times differentiable functions $h$ that*

$$|\mathbb{E}h(W) - \mathbb{E}h(\Sigma^{1/2}Z)| \le \frac{416 d^{7/2}|h|_2 + 960 d^5 |h|_3}{p^{d/2}(1-p)^{3/2}\sqrt{n}}.$$

PROOF. Some rough estimates yield that, for all $1 \le i,j,k \le d$,

$$\lambda^{(i)} \le \frac{15n}{d},$$

$$\operatorname{Var}\mathbb{E}^W(W_i' - W_i)(W_j' - W_j) \le \frac{768 d^5}{n^3 p^d (1-p)^2},$$

$$\mathbb{E}|(W_i' - W_i)(W_j' - W_j)(W_k' - W_k)| \le \frac{64 d^3}{n^{3/2} p^{d/2} (1-p)^{3/2}}.$$

Now apply Theorem 2.1. Details can be found in Appendix B. □

REMARK 4.2. Although the bound is quite crude with respect to the dimension and hence mainly of theoretical interest, it is explicit. For small values of $p$ or large values of $d$, however, Poisson approximation is more appropriate, and in these cases the bounds for normal approximation cannot be expected to be good unless $n$ is very large. We also note that $V_d$ exhibits a local dependence structure and thus also Stein's method using the local approach, such as in Rinott and Rotar (1996), could easily be used; and, of course, there is an abundance of results about $m$-dependent sequences.

REMARK 4.3. In the case of 2-runs, using the notation of (1.8) and the consequent paragraph, it is not difficult to see that, for any choice of $\lambda$ and defining $R = R(V_2, V_1) := \sigma^{-1}(\lambda V_2 - \frac{2}{n}V_2 + \frac{2p}{n}V_1)$, we have that $\lambda^{-1}\sqrt{\operatorname{Var} R}$ is at least of order 1 as $n \to \infty$, where $\sigma^2 := \operatorname{Var} V_2$. It may nevertheless be possible to choose $\lambda$ such that, with $\tilde{R} = \tilde{R}(V_2) := \mathbb{E}^{V_2} R = \sigma^{-1}(\lambda V_2 - \frac{2}{n}V_2 + \frac{2p}{n}\mathbb{E}^{V_2}V_1)$, we have $\lambda\sqrt{\operatorname{Var}\tilde{R}} = \mathrm{o}(1)$, so that a representation (1.2) could indeed be found with $R$ being of the required small order. But, whereas $\mathbb{E}^{V_2}V_1$ is hard to calculate, in this situation the application of the multivariate version (1.7) and Theorem 2.1 is straightforward.

4.3. *Double-indexed permutation statistics.* Let $a_{i,j,k,l}$, $1 \le i,j,k,l \le n$, be real numbers such that $a_{i,j,k,l} = 0$ whenever $i = j$ but $k \ne j$. Assume that

$$\sum_{i,j,k,l} a_{i,j,k,l} = 0 \tag{4.8}$$

and define

$$V_0 = V_0(\pi) = \sum_{s,t=1}^{n} a_{s,t,\pi(s),\pi(t)},$$

where $\pi$ is a uniformly drawn random permutation of size $n$. A Berry–Esseen bound for the distribution of $V_0$ was proved by Zhao et al. (1997) under quite general conditions, generalizing the proof of Bolthausen (1984), which is related to the exchangeable pair coupling. Under similar conditions as Zhao et al. (1997), Barbour and Chen (2005) used the exchangeable pair coupling to find a nontrivial representation of $V_0$ of the form (1.2) with a nonzero remainder term $R$; see their article also for a historical overview. Yet the problem is so rich that there is to date no result which unifies all the cases in which asymptotic normality holds. For example, the results in Barbour and Chen (2005) and in Zhao et al. (1997) do not cover the the number of descents in a random permutation, for which asymptotic normality was derived in Fulman (2004) via exchangeable pairs.

We will discuss here only the applicability of this example to Theorem 2.1 to illustrate the embedding method, which contrasts with Barbour and Chen (2005) in the sense that, with our approach, again one does not need to find a one-dimensional representation of the form (1.2) but can use directly the multidimensional version (1.7) in a straightforward manner. We also do not bound the error terms because the corresponding calculations are too involved for the purpose of this paper.

Construct now an exchangeable pair as follows. Let $I$ and $J$ be distributed uniformly over $1,\ldots,n$ conditioned that $I \ne J$. Define the permutation $\pi' = (\pi(I)\pi(J)) \circ \pi$ so that $\pi'$ is the permutation where $\pi'(k) = \pi(k)$ for all $k \ne I,J$, and where $\pi'(I) = \pi(J)$ and $\pi'(J) = \pi(I)$. Let now, for the sake of a simpler notation, $a^{\pi}_{i,j,k,l} := a_{i,j,\pi(k),\pi(l)}$. Defining $W' = W(\pi')$, we have

$$\begin{aligned} V_0' - V_0 = &- \sum_{s=1}^{n} (a^{\pi}_{I,s,I,s} + a^{\pi}_{J,s,J,s} + a^{\pi}_{s,I,s,I} + a^{\pi}_{s,J,s,J}) \\ &+ (a^{\pi}_{I,I,I,I} + a^{\pi}_{I,J,I,J} + a^{\pi}_{J,I,J,I} + a^{\pi}_{J,J,J,J}) \\ &+ \sum_{s=1}^{n} (a^{\pi}_{I,s,J,s} + a^{\pi}_{J,s,I,s} + a^{\pi}_{s,I,s,J} + a^{\pi}_{s,J,s,I}) \\ &- (a^{\pi}_{I,I,J,J} + a^{\pi}_{I,J,J,I} + a^{\pi}_{J,I,I,J} + a^{\pi}_{J,J,I,I}). \end{aligned}$$

Hence,

$$
\begin{aligned}
&\mathbb{E}^\pi(V_0' - V_0) \\
&\quad = -\frac{1}{n(n-1)} \sum_{i \neq j} \sum_{s=1}^{n} (a^\pi_{i,s,i,s} + a^\pi_{j,s,j,s} + a^\pi_{s,i,s,i} + a^\pi_{s,j,s,j}) \\
&\qquad + \frac{1}{n(n-1)} \sum_{i \neq j} (a^\pi_{i,i,i,i} + a^\pi_{i,j,i,j} + a^\pi_{j,i,j,i} + a^\pi_{j,j,j,j}) \\
&\qquad + \frac{1}{n(n-1)} \sum_{i \neq j} \sum_{s=1}^{n} (a^\pi_{i,s,j,s} + a^\pi_{j,s,i,s} + a^\pi_{s,i,s,j} + a^\pi_{s,j,s,i}) \\
&\qquad - \frac{1}{n(n-1)} \sum_{i \neq j} (a^\pi_{i,i,j,j} + a^\pi_{i,j,j,i} + a^\pi_{j,i,i,j} + a^\pi_{j,j,i,i}) \\
&\quad = -\frac{4}{n} V_0 + \frac{2}{n(n-1)} \sum_{s=1}^{n} \sum_{i \neq j} (a^\pi_{i,s,j,s} + a^\pi_{s,i,s,j}) \\
&\qquad + \frac{2}{n(n-1)} \sum_{i \neq j} (a^\pi_{i,i,i,i} + a^\pi_{i,j,i,j}) - \frac{2}{n(n-1)} \sum_{i \neq j} (a^\pi_{i,i,j,j} + a^\pi_{i,j,j,i}) \\
&\quad = \lambda\left(-\frac{2n-1}{n} V_0 + V_1 + V_2\right) + R_1 + R_2,
\end{aligned}
$$

with $\lambda := 2/(n-1)$ and where

$$
R_1 := \lambda \sum_{i=1}^{n} a^\pi_{i,i,i,i} - \frac{\lambda}{n} \sum_{i,j=1}^{n} a_{i,i,j,j}, \qquad R_2 := -\frac{\lambda}{n} \sum_{i,j=1}^{n} a^\pi_{i,j,j,i},
$$

$$
V_i := \sum_{s=1}^{n} a^{(i)}_{s,\pi(s)} \qquad \text{for } i = 1, 2, \text{ where}
$$

$$
a^{(1)}_{s,t} := \frac{1}{n} \sum_{i,j} a_{s,i,t,j}, a^{(2)}_{s,t} := \frac{1}{n} \sum_{i,j} a_{i,s,j,t}.
$$

Thus, the conditional expectation $\mathbb{E}^\pi(V_0' - V_0)$ can be decomposed into a sum of the original statistic $V_0$ and two related single-indexed permutation statistics, together with an error term. Now, for $i = 1, 2$,

$$
V_i' - V_i = -a^{(i)}_{I,\pi(I)} - a^{(i)}_{J,\pi(J)} + a^{(i)}_{I,\pi(J)} + a^{(i)}_{J,\pi(I)}
$$

and, thus,

$$
\mathbb{E}^\pi(V_i' - V_i) = -\frac{2}{n} V_i + \frac{2}{n(n-1)} \sum_{i \neq j} a^{(i)}_{i,\pi(j)}
$$

$$= -\lambda V_i + \frac{2}{n(n-1)} \sum_{i,j} a^{(i)}_{i,\pi(j)}$$

$$= -\lambda V_i,$$

where the last equality follows from (4.8). Thus, (1.7) holds for the vector $W = (V_0, V_1, V_2)^t$ with

$$\Lambda = \lambda \begin{pmatrix} \frac{2n-1}{n} & -1 & -1 \\ 0 & 1 & 0 \\ 0 & 0 & 1 \end{pmatrix}$$

and $R = (R_1 + R_2, 0, 0)^t$.

In the special case where $a_{ijkl} = b_{ij}c_{kl}$ with $b_{ii} = c_{ii} = 0$ for all $i, j, k, l$ and where $(b_{ij})$ or $(c_{kl})$ is symmetric up to a (possibly negative) constant, we have $R_1 = 0$ and $R_2 = \beta\lambda n^{-1} V_0$ for some number $\beta$, so that (1.7) holds with an $R = 0$ and a slightly different $\Lambda$, which would simplify the estimates. Note that these assumptions hold, for example, if either $(b_{ij})$ or $(c_{ij})$ is the adjacency matrix of an undirected graph containing no self-loops.

*Mann–Whitney–Wilcoxon statistic.* Let $x_1, \ldots, x_{n_x}$ and $y_1, \ldots, y_{n_y}$, $n_x + n_y = n$, be independent random samples from unknown distributions $F_X$ and $F_Y$, respectively. The MWW-statistic is then defined to be the number of pairs $(x_i, y_j)$ such that $x_i < y_j$. Let $\pi(i)$ be the rank of $z_i$, where $z = (x_1, \ldots, x_{n_x}, y_1, \ldots, y_{n_y})$ is the combined sample. To test the hypothesis $H_0 : F_X = F_Y$, we may assume that $\pi$ has uniform distribution. It is easy to see that, defining

$$a_{i,j,k,l} = \begin{cases} +\frac{1}{2}, & \text{if } 1 \le i \le n_x,\ n_x + 1 \le j \le n \text{ and } 1 \le k < l \le n, \\ -\frac{1}{2}, & \text{if } 1 \le i \le n_x,\ n_x + 1 \le j \le n \text{ and } 1 \le l < k \le n, \\ 0, & \text{else,} \end{cases}$$

$V_0$ is equivalent to the MWW-statistic (up to a shift). It is well known that $\operatorname{Var} V_0 = n_x n_y (n+1)/12$ [see Mann and Whitney (1947)], so that if, for some $0 < \alpha < 1$, $n_x \asymp \alpha n$ and $n_y \asymp (1-\alpha)n$, respectively, we have $\operatorname{Var} V_0 \asymp n^3$.

Note now that, as $a_{i,i,k,l} = 0$ for all $i, k, l$ and as $\sum_{i,j} a_{i,j,\pi(j),\pi(i)} = -\sum_{i,j} a_{i,j,\pi(i),\pi(j)}$, we have $R_1 = 0$ and $R_2 = -\frac{\lambda}{n} V_0$. Hence, the remainder term $C$ in Theorem 2.1 has the required lower order.

Further, we calculate that $a^{(1)}_{i,j} = \frac{n_y(n-2j+1)}{2n}$ if $1 \le i \le n_x$ and $a^{(1)}_{i,j} = 0$ otherwise, and therefore, using the variance formula for the usual singly indexed permutation statistics [see Hoeffding (1951)],

$$\operatorname{Var} V_1 = \frac{1}{n-1} \sum_{i,j=1}^{n} (a^{(1)}_{i,j} - a^{(1)}_{i,\cdot} - a^{(1)}_{\cdot,j} + a^{(1)}_{\cdot,\cdot})^2 \asymp n^3.$$

The same asymptotic is true for $V_2$, so that indeed $W = n^{-3/2}(V_0, V_1, V_2)$ with the above coupling and choice of $\Lambda$ is a good candidate for Theorem 2.1.

**5. Some comments on the exchangeability condition.** Exchangeability is used twice in the proof of Theorem 2.1, namely, in (2.7) and (2.10). In this section we discuss the necessity of this condition if one uses the Stein operator of the form in equation (2.4).

5.1. *Exchangeability and anti-symmetric functions.* In (2.7), we use exchangeability in the spirit of Stein (1986). It has been proved by Röllin (2008) that in the one-dimensional setting the exchangeability condition can be omitted for normal approximation by replacing the usual anti-symmetric function (2.6) with $F(w,w') = g(w') - g(w)$, where now only equality in distribution is needed to obtain an identity similar to (2.7). Chatterjee and Meckes (2008) also proved their results with this new function $F$ but under the stronger condition (1.4). However, there seems to be no obvious way to apply the above approach under the more general assumption (1.7) (even with $R = 0$) to remove the exchangeability condition. To see this, note that, by multivariate Taylor expansion,

$$\begin{aligned} g(w') = g(w) + (w' - w)^t \nabla g(w) + \tfrac{1}{2} \nabla^t (w' - w)(w' - w)^t \nabla g(w) \\ + r(w', w), \end{aligned} \tag{5.1}$$

where $r$ is the corresponding remainder term of the expansion. Thus, (5.1) and (1.7) yield the identity

$$\begin{aligned} 0 &= \mathbb{E}g(W') - \mathbb{E}g(W) \\ &= -\mathbb{E}\{W^t \Lambda^t \nabla g(W)\} + \tfrac{1}{2}\mathbb{E}\{\nabla^t (W' - W)(W' - W)^t \nabla g(W)\} \\ &\quad + \mathbb{E}r(W', W), \end{aligned} \tag{5.2}$$

for any suitable function $g$. To optimally match (5.2) and the left-hand side of (2.4), we have to choose $g$ such that the system

$$\Lambda^t \nabla g = \nabla f \tag{5.3}$$

is satisfied. In the one-dimensional setting of Röllin (2008) and the multivariate setting $\Lambda = \lambda \,\mathrm{Id}$ of Chatterjee and Meckes (2008), (5.3) can be solved by setting $g = \lambda^{-1} f$. Indeed, (5.3) cannot be solved in general; only if $\Lambda = \lambda \,\mathrm{Id}$ does (5.3) have a twice continuously partially differentiable solution $g$ for a sufficiently large class of functions $f$.

5.2. *Exchangeability, the covariance matrix and the $\Lambda$ matrix.* In (2.10), using only equality in distribution instead of exchangeability, we obtain

$$\mathbb{E}(W' - W)(W' - W)^t = \Lambda\Sigma + \Sigma\Lambda^t. \tag{5.4}$$

It is clear from (2.11) that the canonical choice for the variance structure of the approximating multivariate normal distribution would then be

$$\tfrac{1}{2}\mathbb{E}(W' - W)(W' - W)^t \Lambda^{-t} = \tfrac{1}{2}(\Lambda\Sigma\Lambda^{-t} + \Sigma) =: \tilde{\Sigma}, \tag{5.5}$$

which in the exchangeable setting reduces to $\Sigma$; see (2.10).

It is easy to see that $\tilde{\Sigma} = \Sigma$ if and only if $\hat{\Lambda} := \Sigma^{-1/2}\Lambda\Sigma^{1/2}$, arising from standardization (see Remark 2.4), is symmetric. If $(W', W)$ is exchangeable, we have from (2.10) that $\tilde{\Sigma} = \Sigma$ and, hence, $\hat{\Lambda}$ is symmetric. While exchangeability of $(W, W)$ is not a necessary condition for $\hat{\Lambda}$ to be symmetric, the following example illustrates that nonsymmetric $\hat{\Lambda}$ is far from unusual.

EXAMPLE 5.1. Let $d$ be a positive integer, $d \geq 4$. Let $X(k) = \{X_i(k); i = 1, \ldots, d\}$, $k \in \mathbb{Z}_+$ be a discrete time Markov chain with values in $\{-1, 1\}^d$ and with the following transition rule. At every time step $k$, pick uniformly an index $I$ from $\{1, 2, \ldots, d\}$. Then with probability $1/2$, let $X_I(k+1) = -X_{I-1}(k)$, and with probability $1/2$, let $X_I(k+1) = X_{I+1}(k)$, where we interpret the indices 0 and $d+1$ as $d$ and 1, respectively. For all $j \neq I$, put $X_j(k+1) = X_j(k)$. Observe that, for arbitrary $k$ and $i \neq j$,

$$\begin{aligned}
&\mathbb{E}[X_i(k+1)X_j(k+1)|X(k)] \\
&\quad = \frac{1}{2d}(X_{i+1}(k) - X_{i-1}(k))X_j(k) + \frac{1}{2d}X_i(k)(X_{j+1}(k) - X_{j-1}(k)) \\
&\qquad + \frac{d-2}{d}X_i(k)X_j(k).
\end{aligned}$$

Now, if $\mathbb{E}\{X_i(k)X_j(k)\} = 0$ for all $i \neq j$, then also $\mathbb{E}\{X_i(k+1)X_j(k+1)\} = 0$ (where the case $j \in \{i-1, i+1\}$ is slightly different than for the other $j$). Thus, if we start the chain such that the $X_i$ are uncorrelated and centered, then, by induction, the $X_i$ are uncorrelated for every $k$ and it is easy to see from this that also the equilibrium distribution of the chain has uncorrelated $X_i$.

Assume that $X^{(1)}, X^{(2)}, \ldots$ is a sequence of mean zero independent and identically distributed $d$-vectors with finite $\Sigma := \mathbb{E}\{X^{(1)}(X^{(1)})^t\}$. It is clear from the multivariate CLT [see, for example, Rotar (1997), page 363, Theorem 4] that $W = n^{-1/2}\sum_{i=1}^n X^{(i)}$ converges to the multivariate mean zero normal distribution with covariance matrix $\Sigma$.

However, consider the following coupling construction. Let $X^{(i)}$ have the equilibrium distribution of the above Markov chain and for each $i$ let $X'^{(i)}$ be the value after one step ahead in the Markov chain, such that the pairs $(X^{(i)}, X'^{(i)})$ are independent for different $i$. Define now $W' = W + n^{-1/2}(X'^{(I)} - X^{(I)})$, where $I$ is uniformly distributed on $\{1, \ldots, n\}$, and note that $\mathscr{L}(W') = \mathscr{L}(W)$. We calculate that $\mathbb{E}^{X^{(i)}}(X'^{(i)} - X^{(i)}) = -\Lambda X^{(i)}$ with

$$\Lambda_{ij} = \frac{1}{d} \cdot \begin{cases} 1, & \text{if } j = i, \\ \frac{1}{2}, & \text{if } j = i-1, \\ -\frac{1}{2}, & \text{if } j = i+1, \\ 0, & \text{else.} \end{cases}$$

Then $\mathbb{E}^W(W' - W) = -n^{-1}\Lambda W$. As $\Lambda$ is not symmetric, $(W, W')$ cannot be exchangeable, and so Theorem 2.1 cannot be applied with this coupling.

## APPENDIX A: PROOF OF COROLLARY 3.1

For $h \in \mathcal{H}$ define the following smoothing:

$$h_s(x) = \int_{\mathbb{R}^d} h(s^{1/2}y + (1-s)^{1/2}x)\Phi(dy), \qquad 0 < s < 1.$$

The following key result for this smoothing can be found in Götze (1991).

Lemma A.1. *Let $Q$ be a probability measure on $\mathbb{R}^d$, and let $W \sim Q, Z \sim \Phi$. Let $a$ be as in (3.1). Then there exists a constant $\gamma > 0$ which depends only on the dimension $d$ such that, for $0 < t < 1$,*

$$\sup_{h\in\mathcal{H}} |\mathbb{E}h(W) - \mathbb{E}h(Z)| \le \gamma\Big[\sup_{h\in\mathcal{H}} |\mathbb{E}(h - \Phi h)_t(W)| + a\sqrt{t}\Big].$$

To prove Corollary 3.1, first we assume that $\Sigma = \mathrm{Id}$. Let $0 < t < 1$. The solution of (2.4) for $h_t$ is $\Psi_t(x) = \frac{1}{2}\int_t^1 \frac{h_s(x) - \Phi h}{1-s}\,ds$, and for $|h| \le 1$, it is shown in Götze (1991) and also in Loh (2008) that there is a constant $\gamma = \gamma(d)$ depending only on the dimension $d$ such that

$$|\Psi_t|_1 \le \gamma, \qquad |\Psi_t|_2 \le \gamma \log(t^{-1}); \tag{A.1}$$

the $\gamma$ is in general not equal to the $\gamma$ in Lemma A.1. Then, as in (2.11),

$$\begin{aligned} |\mathbb{E}h_t(W) - \mathbb{E}h_t(Z)| &= |\mathbb{E}\{\nabla^t\nabla\Psi_t(W) - W^t\nabla\Psi_t(W)\}| \\ &\le \frac{1}{2}\sum_{m,i,j} |(\Lambda^{-1})_{m,i}\mathbb{E}(W_i' - W_i)(W_j' - W_j)(W_k' - W_k)R_{mjk}| \qquad \text{(A.2)} \\ &\quad + \frac{\gamma}{2}\log(t^{-1})A + \gamma C(1 + d\log(t^{-1})), \end{aligned}$$

with $A, B$ and $C$ as in Theorem 2.1. For the last step we used the same estimates as applied for the remainder term in (2.11), and that $\Sigma = \mathrm{Id}$.

For the remainder term $R_{mjk}$, in Loh (2008), Lemma 1 (page 1992), it is shown that, if $|h| \le 1$, then there is a constant $c_0$ (depending only on $d$) such that, for any finite signed measures $Q$ on $\mathbb{R}^d$,

$$\begin{aligned} &\sup_{1\le p,q,r\le d}\left|\int_{\mathbb{R}^d} \frac{\partial^3}{\partial z_p\,\partial z_q\,\partial z_r}\Psi_t(z)Q(dz)\right| \\ &\qquad \le \frac{c_0}{\sqrt{t}} \sup_{0\le s\le 1, y\in\mathbb{R}^d}\left|\int_{\mathbb{R}^d} h(sv + y)Q(dv)\right|. \end{aligned}$$

Thus, we can bound the second term in (A.2) by $\frac{c_0}{2\sqrt{t}}B$. For simplicity, we relabel $\gamma$ as the maximum of $\gamma$, $\gamma^2$ and $\gamma c_0$, yielding that

$$\sup_{h\in\mathcal{H}} |\mathbb{E}h(W) - \mathbb{E}h(Z)| \leq \gamma^2\left(D\log(t^{-1}) + \frac{1}{2}Bt^{-1/2} + C + a\sqrt{t}\right),$$

with $D = \frac{S}{2} + Cd$. The minimum with respect to $t$ is attained for $T = \frac{1}{a^2}(D + \sqrt{\frac{aB}{2} + D^2})^2$, which gives the assertion for $\Sigma = \mathrm{Id}$.

To complete the proof for general $\Sigma$, we standardize

$$Y = \Sigma^{-1/2}W.$$

From condition (C2), we have that for any $d\times d$ matrix $A$ and any vector $b\in\mathbb{R}^d$, $h(Ax+b)\in\mathcal{H}$, so, in particular, $h(\Sigma^{-1/2}x)\in\mathcal{H}$. Hence, the above bounds (A.1) can be applied directly. The proof now continues as for the $\Sigma = \mathrm{Id}$ case, but with the standardized variables. We omit the details.

## APPENDIX B: DETAILS OF THE RUNS EXAMPLE

We first show the following lemma, which may be useful when the nondiagonal entries of $\Lambda$ are small compared to the diagonal-entries.

LEMMA B.1. *Assume that $\Lambda$ is lower triangular and assume that there is $a > 0$ such that $|\Lambda_{i,j}| \leq a$ for all $j < i$. Then, with $\gamma := \inf_i |\Lambda_{ii}|$,*

$$\sup_i \lambda^{(i)} \leq \frac{(a/\gamma + 1)^{d-1}}{\gamma}.$$

PROOF. Write $\Lambda = \Lambda_E\Lambda_D$, where $\Lambda_D$ is diagonal with the same diagonal as $\Lambda$ and $\Lambda_E$ is lower triangular with diagonal entries equal to 1 and $(\Lambda_E)_{i,j} := \Lambda_{i,j}/\Lambda_{j,j}$. Denote by $\|\cdot\|_p$ the usual $p$-norm for matrices and recall that, for any matrix $A$, $\|A\|_1 = \sup_j \sum_i |A_{i,j}|$. Then, $\lambda^{(i)} \leq \|\Lambda^{-1}\|_1 \leq \|\Lambda_D^{-1}\|_1\|\Lambda_E^{-1}\|_1$. Noting that $|(\Lambda_E)_{i,j}| \leq a/\gamma$ for all $j < i$, we have from Lemeire (1975) that $\|\Lambda_E^{-1}\|_1 \leq (a/\gamma+1)^{d-1}$. Now, as $\|\Lambda_D^{-1}\|_1 = \gamma^{-1}$, the claim follows. □

Fix now $i$ and $j$. From (4.2) it is not difficult to see that we can find two sequences $A_1,\ldots,A_{N_{i,j}}$ and $B_1,\ldots,B_{N_{i,j}}$ of subsets of $\{-d+1,\ldots,2d-3\}$ such that

$$\begin{aligned}\mathbb{E}^{\xi,\xi'}(V_i'-V_i)(V_j'-V_j) &= \frac{1}{n}\sum_{m=1}^{n}\sum_{k=1}^{N_{i,j}}\prod_{l\in A_k}\xi_{m+l}\prod_{l\in B_k}\xi'_{m+l}\\ &=: \frac{1}{n}\sum_{m=1}^{n}\nu^{i,j}(m).\end{aligned}\tag{B.1}$$

From (4.2) is easy to see that $N_{i,j} \le 4(d+i-2)(d+j-2) \le 16d^2$, as $V_i' - V_i$ (respectively $V_j' - V_j$) contain no more than $2(d+i-2)$ [respectively $2(d+j-2)$] summands. Note that $|A_k| + |B_k| \ge i \vee j$, that is, every summand in (B.1) is the product of at least $i \vee j$ independent random indicators. Hence, it is not difficult to see that

$$\operatorname{Var}(\nu^{i,j}(m)) \le 256 d^4 p^{i\vee j}. \tag{B.2}$$

Now,

$$\begin{aligned}
&\operatorname{Var}\mathbb{E}^W (W_i' - W_i)(W_j' - W_j)\\
&\quad\le \frac{1}{n^2 p^{i+j}(1-p)^2} \operatorname{Var}\mathbb{E}^{\xi,\xi'}(V_i' - V_i)(V_j' - V_j)\\
&\quad= \frac{1}{n^4 p^{i+j}(1-p)^2} \sum_{m,m'=1}^{n} \operatorname{Cov}(\nu^{i,j}(m), \nu^{i,j}(m')).
\end{aligned}$$

If $|m - m'| \ge 3d$, we have $\operatorname{Cov}(\nu^{i,j}(m), \nu^{i,j}(m')) = 0$ because $\nu^{i,j}(m)$ and $\nu^{i,j}(m')$ are independent. If $|m - m'| < 3d$, we can apply (B.2) to estimate the covariances and, hence, we obtain

$$\operatorname{Var}\mathbb{E}^W (W_i' - W_i)(W_j' - W_j) \le \frac{768 d^5}{n^3 p^{i\wedge j}(1-p)^2}.$$

Similar arguments lead to the estimate

$$\mathbb{E}|(V_i' - V_i)(V_j' - V_j)(V_k' - V_k)| \le 64 d^3 p^{\max\{i,j,k\}},$$

hence, for the second summand in (2.2),

$$\mathbb{E}|(W_i' - W_i)(W_j' - W_j)(W_k' - W_k)| \le \frac{64 d^3 p^{\max\{i,j,k\}}}{n^{3/2} p^{(i+j+k)/2}(1-p)^{3/2}}.$$

Applying Lemma B.1 to the matrix $n\Lambda$ with $a = 2$ and $\gamma = d - 1$, we obtain

$$\lambda^{(i)} \le \frac{n(2/(d-1)+1)^{d-1}}{(d-1)} \le \frac{15n}{d}.$$

Combining all estimates with Theorem 2.1 proves Theorem 4.1.

**Acknowledgments.** The authors would like to thank two anonymous referees for their helpful comments, which led to an improvement of the paper. A. R. would like to thank the Department of Statistics at the University of Oxford for the kind support during the time when most of the research for this paper was done.

Department of Statistics
University of Oxford
1 South Parks Road
Oxford OX1 3TG
United Kingdom
E-mail: reinert@stats.ox.ac.uk

Department of Statistics
National University of Singapore
2 Science Drive 2
Singapore 117543
E-mail: staar@nus.edu.sg